\documentclass[12pt]{amsart}
\usepackage{amssymb, colordvi}
\usepackage{multirow}
\usepackage{graphicx}

\numberwithin{equation}{section}
\newtheorem{thm}{Theorem}
\numberwithin{thm}{section}

\newtheorem{lemma}[thm]{Lemma}
\newtheorem{cor}[thm]{Corollary}

\newtheorem{example}[thm]{Example}
\newtheorem{remark}[thm]{Remark}

\newenvironment{rem}{\begin{remark}\rm}{\end{remark}}
\newcounter{FNC}[page]
\def\fauxfootnote#1{{\addtocounter{FNC}{2}$^\fnsymbol{FNC}$%
     \let\thefootnote\relax\footnotetext{$^\fnsymbol{FNC}$#1}}}

\copyrightinfo{}{}
\newcommand{\calA}{\mathcal{A}}

\newcommand{\calI}{\mathcal{I}}

\newcommand{\calV}{\mathcal{V}}
\newcommand{\calW}{\mathcal{W}}

\newcommand{\C}{\mathbb{C}}
\newcommand{\N}{\mathbb{N}}

\newcommand{\R}{\mathbb{R}}
\newcommand{\Z}{\mathbb{Z}}
\renewcommand{\P}{\mathbb{P}}

\newcommand{\New}{{\rm New}}
\newcommand{\vol}{{\rm vol}}

\title[Fewnomial Upper Bounds from Gale dual polynomial systems]{New Fewnomial
  Upper Bounds from\\ Gale dual polynomial systems} 

\author{Fr\'ed\'eric Bihan}
\address{Laboratoire de Math\'ematiques\\
         Universit\'e de Savoie\\
         73376 Le Bourget-du-Lac Cedex\\
         France}

\email{Frederic.Bihan@univ-savoie.fr}
\urladdr{http://www.lama.univ-savoie.fr/\~{}bihan/}

\author{Frank Sottile}
\address{Department of Mathematics\\
         Texas A\&M University\\
         College Station\\
         Texas \ 77843\\
         USA}
\email{sottile@math.tamu.edu}
\urladdr{http://www.math.tamu.edu/\~{}sottile/}

\thanks{Sottile supported by the Institut Henri Poincar\'e, NSF CAREER grant
  DMS-0538734, and Peter Gritzmann of the Technische
         Universit\"at M\"unchen}  

\begin{document}

\begin{abstract}
 We show that there are fewer than 
 $\frac{e^2+3}{4}\, 2^{\binom{k}{2}}n^k$ positive solutions to a fewnomial
 system consisting of $n$ polynomials in $n$ variables having a total of
 $n{+}k{+}1$ distinct monomials.
 This is significantly smaller than Khovanskii's fewnomial bound of
 $2^{\binom{n{+}k}{2}}(n{+}1)^{n{+}k}$.
 We reduce the original system to a system of $k$ equations in $k$
 variables which depends upon the vector configuration Gale dual to the
 exponents of the monomials in the original system.
 We then bound the number of solutions to this Gale system.
 We adapt these methods to show that a hypersurface in the positive orthant
 of $\R^n$ defined by a polynomial with $n{+}k{+}1$ monomials has at most
 $C(k)n^{k-1}$ compact connected components.
 Our results hold for polynomials with real exponents.
\end{abstract}
\dedicatory{We dedicate this paper to Askold Khovanskii on the occasion of
his 60th birthday.}
\maketitle

%
\section*{Introduction}
%

The zeroes of a (Laurent) polynomial $f$ in $n$ variables lie in the complex torus
$(\C^\times)^n$ and the support of $f$
is the set $\calW\subset\Z^n$ of exponent vectors of its monomials.
Kouchnirenko~\cite{Ko75} showed that the number of non-degenerate solutions to a
system of $n$ polynomials with support $\calW$ is at most the volume of the
convex hull of $\calW$, suitably normalized.
This bound is sharp for every $\calW$.

It is also a bound for the number of real solutions to real polynomial
equations, but it is not sharp. 
In 1980, Khovanskii~\cite{Kh80} gave a bound for the number of non-degenerate
solutions in the positive orthant $\R_>^n$ to a system of polynomials which
depends only on the cardinality $|\calW|$ of $\calW$.
Multiplying this bound by the number $2^n$ of orthants in $(\R^\times)^n$ gives
a bound for the total number of real solutions which is 
smaller than the Kouchnirenko bound when $|\calW|$ is small
relative to the volume of the convex hull of $\calW$.

Khovanskii's fewnomial bound is quite large.
If $|\calW|=n{+}k{+}1$, then it is
\[
   2^{\binom{n+k}{2}}(n+1)^{n+k}\,.
\]
When $n=k=2$, this becomes $2^6\cdot 3^4 = 5184$.
The first concrete result showing that Khovanskii's bound is likely overstated
was due to Li, Rojas, and Wang~\cite{LRW03} who showed that two trinomials in
two variables have at most 5 positive solutions.
Since we may multiply the polynomials by monomials without changing their
positive solutions, we can assume that they both have a constant term.
Thus there are at most 5 monomials between them and so this is a fewnomial system
with $n=k=2$.
While $5$ is considerably smaller than $5184$,
two trinomials do not constitute a general fewnomial system with
$n=k=2$, and the problem of finding a bound tighter than $5184$ in this case 
remained open.

Our main result is a new fewnomial upper bound, which is a special case of
Theorem~\ref{FinalBounds}.\medskip

\noindent{\bf Theorem.}
{\em 
 A system of $n$ polynomials in $n$ variables having a total of $n{+}k{+}1$
 distinct monomials has fewer than 
 \[
     \tfrac{e^2+3}{4}\, 2^{\binom{k}{2}}n^k \leqno{(*)} 
 \]
 non-degenerate solutions in the positive orthant.}\medskip

When $n=k=2$, this bound is $20$.
Like Khovanskii's fewnomial bound, the bound $(*)$ is 
for polynomials whose exponents in $\calW$ may be real vectors.
We will use real exponents to simplify our proof in a key way.

While this bound $(*)$ is significantly smaller than
Khovanskii's bound, they  both have the same quadratic factor of $k$ in the
exponent of 2. 
Despite that similarity, 
our method of proof is different from Khovanskii's induction.
A construction based on~\cite{Bihan} gives a lower bound of $(1+n/k)^k$ real
solutions.
Thus the dependence, $n^k$, on $n$ in the bound~$(*)$ cannot be
improved.
We believe that the term $2^{\binom{k}{2}}$ is considerably overstated.
Consequently, we do not consider it important to obtain a smaller
constant than $\frac{e^2+3}{4}$, even though there is room in our proof 
to improve this constant.

This constant comes from estimates 
for a smaller bound 
that we prove using polyhedral combinatorics
and some topology of toric varieties.
When $k=2$, this smaller bound becomes 
$2n^2+\bigl\lfloor\frac{(n+3)(n+1)}{2}\bigr\rfloor$,
and when $n=2$ it is 15.
While $15$ is smaller than $5184$, we believe that the true bound is still
smaller. 
\smallskip

When $k{=}1$, 
Bihan~\cite{Bihan} gave the tight bound of $n{+}1$.
This improved an earlier bound of $n^2{+}n$ by Li, Rojas, and Wang~\cite{LRW03}. 
Bihan's work built upon previous work of Bertrand, Bihan, and Sottile~\cite{BBS}
who gave a bound of $2n{+}1$ for all solutions in $(\R^\times)^n$ for Laurent
polynomials with primitive support.
For these results, a general polynomial system supported on $\calW$ is reduced
to a peculiar univariate polynomial, which was analyzed to obtain the 
bounds. 

Our first step is to generalize that reduction.
A system of $n$ polynomials in $n$ variables having a total of
$n{+}k{+}1$ distinct monomials has an associated Gale system,
which is a system of $k$ equations in variables 
$y=(y_1,\dotsc,y_k)$ having the form
 \[
   \prod_{i=1}^{n+k} p_i(y)^{a_{i,j}}\ =\ 1
   \qquad\mbox{for}\quad j=1,2,\dotsc,k\,,  
 \]
where each polynomial $p_i(y)$ is linear in $y\in\R^k_>$
and the exponents $a_{i,j}\in\R$ come from a vector configuration Gale dual
to $\calW$. 
Our bounds are obtained by applying the Khovanskii-Rolle
Theorem~\cite[pp.~42--51]{Kh91} (following some ideas of~\cite{GNS})
and some 
combinatorics of the polyhedron defined by $p_i(y)>0$, for
$i=1,\dotsc,n{+}k$.\smallskip

We also use Gale systems to show that the number of
compact connected components of a hypersurface in $\R^n_>$
given by a polynomial with $n{+}k{+}1$ monomials
is at most
\[
   \left(\tfrac{k}{2}2^{\binom{k}{2}}\,+\,\tfrac{e^2+1}{8}k 2^{\binom{k-1}{2}}
      \right) n^{k-1} 
      \ +\ 
    \left( \tfrac{e^2}{8}2^{\binom{k-2}{2}} \right) n^{k-2}\,,
\]
when $k\leq n$.
When $k=2$, we improve this to $\lfloor\frac{5n+1}{2}\rfloor$.
These results improve earlier bounds.
Khovanskii~\cite[\S 3.14]{Kh80} gives a bound having the order
$2^{\binom{n+k}{2}}(2n^3)^{n-1}$ for the total Betti number
of such a hypersurface, and Perrucci~\cite{Pe05} gave the bound
$2^{\binom{n+k}{2}}(n{+}1)^{n+k}$ for the number of all components of such a
hypersurface. 
\medskip

We first recall some basics on polynomial systems in Section 1.
Section 2 shows how to reduce a polynomial system 
to its Gale dual system.
In Section 3, we establish bounds for the number of positive solutions to a
system with $n{+}k{+}1$ monomials.
Good bounds are known when $n=1$~\cite{D1636} or $k=1$~\cite{Bihan}, so we
restrict attention to when $k,n\geq 2$.
In Section 4 we give bounds on the
number of compact components of a hypersurface in $\R^n_>$ defined by
a polynomial with $n{+}k{+}1$ monomials.

%
\section{Polynomial systems and vector configurations}
%

Fix positive integers $n$ and $k$.
Let $\calW=\{w_0,w_1,\dotsc,w_{n+k}\}\subset\Z^n$ be a collection of integer
vectors, which are exponents for monomials in $z\in(\C^\times)^n$.
A linear combination of monomials with exponents from $\calW$,
 \begin{equation}\label{Eq:sparse_poly}
    f(z)\ =\ c_0z^{w_0}\;+\;c_1z^{w_1}\;+\;
    c_2z^{w_2}\;+\;\dotsb \;+\;  c_{n+k}z^{w_{n+k}}\,,
 \end{equation}
is a (Laurent) polynomial with \Blue{{\it support $\calW$}}.
We will always assume that the coefficients $c_i$
are real numbers.
A \Blue{{\it system with support $\calW$}} is a system 
  \begin{equation}\label{Eq:sparse_system}
   f_1(z)\ =\ f_2(z)\ =\ \dotsb\ =\ f_n(z)\ =\ 0\,,
 \end{equation}
of polynomials, each with support $\calW$.
Solutions to~\eqref{Eq:sparse_system} are points in $(\C^\times)^n$.
Multiplying each polynomial $f_i(z)$ by the monomial $z^{-w_0}$ does not change
the solutions to the system~\eqref{Eq:sparse_system}, but it 
translates $\calW$ by $-w_0$ so that $w_0$ becomes the origin.
We shall henceforth assume that $0\in\calW$ and write
$\calW=\{0,w_1,\dotsc,w_{n+k}\}$.

Kouchnirenko~\cite{Ko75} gave a bound for the number of solutions to a polynomial
system~\eqref{Eq:sparse_system}.
Let $\vol(\calW)$ be the Euclidean volume of the convex hull of $\calW$.\medskip

\noindent{\bf Kouchnirenko's Theorem.}\ 
{\em 
   The number of isolated solutions in $(\C^\times)^n$ to a
   system~$\eqref{Eq:sparse_system}$ with support $\calW$
   is at most $n!\vol(\calW)$.
   This bound is attained for generic systems with support $\calW$.
}\medskip

Here, generic means that there is a non-empty Zariski open subset consisting of
systems where the maximum is acheived.
Every solution $z$ to a generic system is \Blue{{\it non-degenerate}} in that
the differentials $df_1,\dotsc,df_n$ are linearly independent at $z$.
\smallskip

For $w\in\R^n$, $z\mapsto z^w$ is a well-defined analytic function on the
positive orthant $\R^n_>$.
A \Blue{{\it polynomial} ({\it with real exponents})} is a linear combination of such
functions  and its \Blue{{\it support} $\calW}\subset\R^n$ is the set of exponent
vectors.
While Kouchnirenko's Theorem does not apply to systems of polynomials with 
real exponents, Khovanskii's bound does.
Our bounds are likewise for polynomials with real exponents.

The usual theory of discriminants for ordinary polynomial systems holds in this
setting (see~\cite[\S 2]{Na01} for more details).
In the space $\R^{2n(n+k+1)}$ of parameters (coefficients and exponents) for
systems of polynomials with real exponents having $n{+}k{+}1$ monomials
there is a discriminant hypersurface $\Sigma$ which is defined by the vanishing
of an explicit analytic function, and has the following property:
In every connected component of the complement of $\Sigma$ the number of
positive solutions is constant, and all solutions are non-degenerate.
Perturbing the coefficients and exponents of such a system does not 
change the number of non-degenerate positive solutions,
so we may assume that the exponents are rational numbers.
After a monomial change of coordinates, we may assume that the exponents are
integers.
Thus the maximum number of non-degenerate positive solutions to
a system of $n$ polynomials in $n$ variables having a total of $n{+}k{+}1$
monomials with real exponents is achieved by Laurent polynomials.\smallskip

We sometimes write $V(f)$ for the set of zeroes of a function $f$ on a given
domain $\Delta$, and $V(f_1,\dotsc,f_j)$ for the common zeroes of functions
$f_1,\dotsc,f_j$ on $\Delta$. 

%
\section{Gale systems}\label{Galesystems}
%

Given a system of $n$ polynomials in $n$ variables with common support
$\calW$ where $|\calW|=n{+}k{+}1$, we use the linear
relations on $\calW$ (its Gale dual configuration)
to obtain a special system of $k$ equations in $k$ variables, called a Gale 
system. 
Here, we treat the case when $\calW$ is a configuration of real vectors,
and prove an equivalence between positive solutions of the original system
and solutions of the associated Gale system. 
There is another version for non-zero complex solutions which holds when $\calW$
is a configuration of integer vectors.
We will present that elsewhere~\cite{BS06}.

Suppose that $\calW=\{0,w_1,w_2,\dotsc,w_{n+k}\}$ spans $\R^n$,
for otherwise no system with support $\calW$ has any non-degenerate solutions.
We assume that the last $n$ vectors are linearly independent.
We consider a system with support $\calW$ such that the coefficient
matrix of $z^{w_1},\dotsc,z^{w_n}$ is invertible, and then put it into
diagonal form
 \begin{equation}\label{Eq:diagonal}
    z^{w_i}\ =\ p_i(z^{w_{n+1}},z^{w_{n+2}},\dotsc,z^{w_{n+k}})\ =:\ g_i(z)
     \qquad\mbox{for}\quad i=1,\dotsc,n\,,
 \end{equation}
where each $p_i$ is linear, so that $g_i(z)$ is a polynomial with support in
$\{0, w_{n+1},\dotsc,w_{n+k}\}$. 

A vector $a\in\R^{n+k}$ with
\[
   0\ =\ w_1 a_1 + w_2 a_2 + \dotsb + w_{n+k} a_{n+k}
\]
gives a monomial identity
\[
   1 \ =\  (z^{w_1})^{a_1} (z^{w_2})^{a_2}\dotsb (z^{w_{n+k}})^{a_{n+k}}\,.
\]
If we substitute the polynomials $g_i(z)$ of~\eqref{Eq:diagonal} into this
identity, we obtain a consequence of the
system~\eqref{Eq:diagonal} ,
 \begin{equation}\label{Eq:Gale_monomial}
   1 \ =\  (g_1(z))^{a_1} \dotsb (g_n(z))^{a_n} \cdot
        (z^{w_{n+1}})^{a_{n+1}} \dotsb (z^{w_{n+k}})^{a_{n+k}}\,.
 \end{equation}

Let $A=(a_{i,j})\in\R^{(n+k)\times k}$ be a matrix whose columns 
form a basis for the linear relations among the elements
of $\calW$.
%
%
Under the substitution
 \begin{equation}\label{Eq:y_subs}
   y_j\ =\ z^{w_{n+j}}\qquad\mbox{for}\quad j=1,\dotsc,k\,,
 \end{equation}
the polynomials $g_i(z)$ become linear functions $p_i(y)$.
Set $p_{n+i}(y):=y_i$ for $i=1,\dotsc,k$.
Then the $k$ equations of the form~\eqref{Eq:Gale_monomial} given by the
columns of $A$
 \begin{equation}\label{Eq:GaleSystem}
   1\ =\ \prod_{i=1}^{n{+}k} p_i(y)^{a_{i,j}}
   \qquad\mbox{for}\quad j=1,\dotsc,k\,,
 \end{equation}
constitute a \Blue{{\it Gale system}} associated to the diagonal
system~\eqref{Eq:diagonal}. 
This system only makes sense for variables $y$ in the
polyhedron
\[
   \Delta\ :=\ \{ y\mid p_i(y)>0\ \mbox{for}\ i=1,\dotsc,n{+}k\}\ 
    \subset \R^k_>\,.
\]

\begin{rem}\label{R:lin=mult}
While the system~\eqref{Eq:GaleSystem} appears to depend upon a choice
of basis (the columns of $A$) for the linear relations among $\calW$, in fact it
does not 
as the invertible linear transformation between two bases gives 
an invertible multiplicative transformation between their corresponding Gale systems.
\end{rem}

\begin{thm}\label{Th:GaleSystem}
 Set $\calV:=\{w_{n+1}, \dotsc,w_{n+k}\}$.  
 Then the association 
\[
  \varphi_\calV\ \colon\ 
   \R^n_>\ \ni\ z\ \longmapsto\ (z^{w_{n+1}},z^{w_{n+2}},\dotsc,z^{w_{n+k}})
   \ =:\ y\ \in\ \R^k_>
\]
 is a bijection between solutions $z\in\R^n_>$ to the diagonal 
 system~\eqref{Eq:diagonal} and solutions $y\in\Delta$ to the
 Gale system~\eqref{Eq:GaleSystem} which restricts to a bijection between
 their non-degenerate solutions.
\end{thm}
%
%

\begin{proof}
 Since the vectors $\{w_{k+1},\dotsc,w_{n+k}\}$ are linearly
 independent, we may change coordinates and assume that
 for $1\leq i\leq n$ we have $w_{k+i}=e_i$, the $i$th standard basis
 vector. 
 The linear relations among the vectors in $\calW$ then have a basis of the form 
 \begin{equation}\label{Eq:LR}
  w_j\ =\ w_{j,1}e_1 +  w_{j,2}e_2 + \dotsb + w_{j,n}e_n
  \qquad\mbox{for}\quad j=1,2,\dotsc,k\,.
 \end{equation}

 Suppose first that $k>n$.
 Then the original system is
\begin{equation}\label{Eq:OS}
   z^{w_j}\ =\ p_j(z^{w_{n+1}}, z^{w_{n+2}}, \dotsc, z^{w_{k}},\ 
                   z_1,\dotsc,z_n) 
                     \qquad\mbox{for}\quad j=1,2,\dotsc,n\,.
 \end{equation}
 Using the linear relations~\eqref{Eq:LR}, the Gale system becomes
 \begin{eqnarray}
   1&=& p_j(y)^{-1} \prod_{i=1}^n y_{k-n+i}^{w_{j,i}}
      \label{Eq:G_orig}
           \qquad\mbox{for}\quad j=1,2,\dotsc,n\,,\qquad\mbox{and}\\
   1&=& y_{j-n}^{-1}\prod_{i=1}^n y_{k-n+i}^{w_{j,i}}
      \label{Eq:G_zees}
           \qquad\mbox{for}\quad j=n{+}1,n{+}2,\dotsc,k\,.
 \end{eqnarray}

%
%
 The image $\varphi_\calV(\R^n_>)\subset\R^k_>$ is a
 graph defined by the subsystem~\eqref{Eq:G_zees}.
 Given this, the subsystem~\eqref{Eq:G_orig} is just a restatement of the
 original system~\eqref{Eq:OS}.
 This proves the theorem in the case $k>n$ as this equivalence of systems
 respects the multiplicity of solutions.\smallskip

 Suppose now that $k\leq n$.
 Then $y_i=z_{n{+}i}$ for $i=1,\dotsc,k$, and~\eqref{Eq:OS} becomes
 \begin{eqnarray}
  z^{w_j}&=& p_j(z_{n-k+1},z_{n-k+2},\dotsc,z_n)
       \label{Eq:Orig_G}
           \qquad\mbox{for}\quad j=1,2,\dotsc,k\,,\\
  z_j&=& p_{j+k}(z_{n-k+1},z_{n-k+2},\dotsc,z_n)
       \label{Eq:Orig_zees}
           \qquad\mbox{for}\quad j=1,2,\dotsc,n{-}k\,.
 \end{eqnarray}
 Using the linear relations~\eqref{Eq:LR}, the Gale system is
 \begin{equation}\label{Eq:N_G}
   1\ =\  p_j(y)^{-1} \cdot
     \prod_{i=1}^{n-k} p_{j+k}(y)^{w_{j,i}}\cdot
     \prod_{i=n-k+1}^n y_{k-n+i}^{w_{j,i}} 
           \qquad\mbox{for}\quad j=1,2,\dotsc,k\,. 
 \end{equation}

 Conversely, the Gale system~\eqref{Eq:N_G} together with the
 subsystem~\eqref{Eq:Orig_zees} implies the subsystem~\eqref{Eq:Orig_G}.
 This proves the theorem when $k\leq n$.
\end{proof}

%
\section{Upper bounds for Gale  systems}\label{S:general_k}
%

Let $n,k\geq 2$ be positive integers and $y_1,\dotsc,y_k$ be indeterminates.
Let $B=(b_{i,j})$ be a $(n{+}k)\times (k{+}1)$ matrix of real numbers whose 
columns are indexed by $0,1,\dotsc,k$.
Define linear polynomials $p_i(y)$ by
\[
   p_i(y)\ :=\ 
    b_{i,0} + b_{i,1}y_1 + \dotsb + b_{i,k}y_k
   \qquad\mbox{for}\quad i=1,2,\dotsc,n{+}k\,,
\]
and define the (possibly unbounded) polyhedron
\[
  \Delta\ :=\ \{y\in\R^k\mid p_i(y)>0\ \mbox{for}\ i=1,2,\dotsc,n{+}k\}\,.
\]

Suppose that $A=(a_{i,j})$ is a $(n{+}k)\times k$ matrix of real numbers, none of
whose row vectors is zero.
We consider the system of equations on $\Delta$ of the form
 \begin{equation}\label{Eq:New_system}
    1\ =\ \Blue{f_j(y)}\ :=\ \prod_{i=1}^{n{+}k} p_i(y)^{a_{i,j}}\,
    ,\qquad\mbox{for}\quad j=1,2,\dotsc,k\,.
 \end{equation}
%

\begin{thm}\label{Bounds}
  The system~\eqref{Eq:New_system} has fewer than 
 $\tfrac{e^2+3}{4} \, 2^{\binom{k}{2}}n^k$ non-degenerate solutions in $\Delta$.
\end{thm}

Given a set 
$\calW=\{0,w_1,\dotsc,w_{n+k}\} \subset \R^n$ wich spans $\R^n$,
let $N$ be the number of non-zero rows in a
$(n{+}k)\times k$ matrix of linear relations among 
$\{w_1,\dotsc,w_{n+k}\}$ and set
$n(\calW):=N-k$.
Then $n(\calW)\leq n$.

\begin{thm}\label{FinalBounds}
 A polynomial system supported on 
 $\calW=\{0,w_1,\dotsc,w_{n+k}\} \subset \R^n$ has fewer than
 \[
    \tfrac{e^2+3}{4} \, 2^{\binom{k}{2}} n(\calW)^k
 \]
 non-degenerate positive solutions.
\end{thm}

\begin{proof}[Proof of Theorem~{$\ref{FinalBounds}$}]
 Suppose that the vectors $\calW$ are ordered so that $w_{k+1},\dotsc,w_{n+k}$
 are linearly independent.
 Given a system supported on $\calW$, we perturb its coefficients 
 without decreasing its number of non-degenerate solutions, 
 put it into diagonal form~\eqref{Eq:diagonal}, and then consider an 
 associated Gale system.
 By Theorem~\ref{Th:GaleSystem}, the number of non-degenerate positive 
 solutions to the diagonal system is equal to the number of non-degenerate 
 solutions to the Gale system in the polyhedron $\Delta$. 
 This Gale system has the form~\eqref{Eq:New_system} where the number of 
 factors is $n(\calW)+k$, and so the theorem follows from 
 Theorem~\ref{Bounds}.
\end{proof}

We will want to consider systems that are equivalent
to~\eqref{Eq:New_system}. 
By Remark~\ref{R:lin=mult}, we obtain an equivalent system from any matrix 
$A'$ of size 
$(n{+}k)\times k$ with the same column span as $A$.
Since we assumed that no row vector of $A$ vanishes,
we may thus assume that no entry of the matrix vanishes, and
hence no exponent $a_{i,j}$ vanishes.

We may perturb the (real-number) exponents $a_{i,j}$ in~\eqref{Eq:New_system}
without decreasing the number of non-degenerate solutions. 
Doing this, if necessary, we may assume that every square submatrix of $A$ is
invertible, and that the exponents are all rational numbers.
We will invoke this assumption to insure that solution sets to
certain equations are smooth algebraic varieties and have the expected dimension.
Similarly, we may perturb the coefficients $B$ of the linear factors $p_i(y)$ 
so that they are in general position
and thus the polyhedron $\Delta$ is \Blue{{\it simple}} in that every face of 
codimension $j$ lies on exactly $j$ facets.

To prove Theorem~\ref{Bounds} we use Khovanskii's generalization
of Rolle's Theorem to 1-manifolds (see \S3.4 in~\cite{Kh91}).
Let $f_1,\dotsc,f_k$ be smooth functions defined on $\Delta\subset\R^k$
which have finitely many common zeroes $V(f_1,\dotsc,f_k)$ in $\Delta$, 
where $V(f_1,\dotsc,f_{k-1})$ is a smooth curve $C$ in $\Delta$.
Let $\flat(C)$ denote the number of unbounded components of $C$ in $\Delta$, and
let
\[
   \Gamma\ =\ J(f_1,f_2,\dotsc,f_k)\ :=\ 
   \det\left( 
         \begin{matrix} \frac{\partial f_j}{\partial y_l}\end{matrix} 
           \right)_{j,l=1,\dotsc,k}
\]
be the Jacobian of $f_1,\dotsc,f_k$.
This vanishes when the differentials $df_1(y),\dotsc,df_k(y)$
are linearly dependent.

\begin{thm}[Khovanskii-Rolle]\label{T:Kh-Ro}
  $ |V(f_1,f_2,\dotsc,f_k)| \ \leq\ \flat(C) + |V(f_1,f_2,\dotsc,f_{k-1}, \Gamma)|$. 
\end{thm} 

The simple idea behind this basic estimate is that along any arc of the curve
$C$ connecting two points, say $a$ and $b$, where $C$ meets the hypersurface
$\{y\mid f_k(y)=0\}$, there must be a point $c\in C$ between $a$ and $b$ where
$C$ is tangent to a level set of $f_k$, so that the normal $df_k$ to this level
set is normal to $C$.
Since the normal directions to $C$ are spanned by $df_1,\dotsc,df_{k-1}$, this
means that $\Gamma$ vanishes at such a point $c$.

\begin{figure}[htb]
\[
  \begin{picture}(270,130)
   \put(0,0){\includegraphics{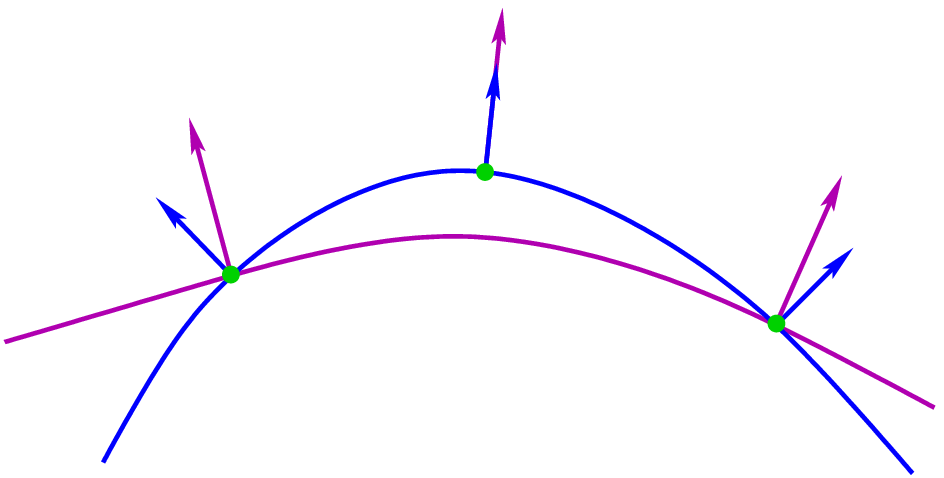}}
   \put(40,9){\Blue{$C=V(f_1)$}} \put( 2,28){\Magenta{$V(f_2)$}}
   \put(65,48){\ForestGreen{$a$}}
   \put(31,69){\Blue{$df_1$}}  \put(60,97){\Magenta{$df_2$}}
   \put(136,77){\ForestGreen{$c$}}
   \put(122,111){\Blue{$df_1$}}  \put(149,124){\Magenta{$df_2$}}
   \put(217,33){\ForestGreen{$b$}}
  \put(245,55){\Blue{$df_1$}}  \put(222,79){\Magenta{$df_2$}}
  \end{picture}
\]
\caption{$\Magenta{df_1}\wedge\Blue{df_2}(\ForestGreen{a})>0$,
         $\Magenta{df_1}\wedge\Blue{df_2}(\ForestGreen{c})=0$, and 
         $\Magenta{df_1}\wedge\Blue{df_2}(\ForestGreen{b})<0$.}
\label{F:KhRoEx}
\end{figure}

We can use the ordinary Rolle theorem to prove Theorem~\ref{T:Kh-Ro}.
Suppose
that the arc along $C$ from $a$ to $b$ is parametrized by a
differentiable function 
$\varphi\colon [0,1]\to C$ with $\varphi(0)=a$ and $\varphi(1)=b$.
By Rolle's Theorem applied to the function $f(t):=f_k(\varphi(t))$, we obtain a
point $c\in(a,b)$ such that $f'(c)=0$, and so $\Gamma$ vanishes at $\varphi(c)$.
Figure~\ref{F:KhRoEx} illustrates the Khovanskii-Rolle Theorem when $k=2$.

On $\Delta$, the system~\eqref{Eq:New_system} has an equivalent formulation in
terms of logarithmic equations
 \begin{equation}\label{Eq:Lop_equation}
  0\ =\ \psi_j(y)\ :=\ 
      \log\bigl(\prod_{i=1}^{n+k} p_i(y)^{a_{i,j}}\bigr)
     \ =\  {\displaystyle\sum_{i=1}^{n+k}} a_{i,j}\log (p_i(y))\qquad 
    \mbox{for}\quad j=1,2,\dotsc,k\,.
 \end{equation}
Define functions $\Gamma_k(y),\dotsc,\Gamma_1(y)$ by the recursion 
\[
   \Gamma_j\ :=\ J(\psi_1,\psi_2,\dotsc,\psi_j,\, 
              \Gamma_{j+1},\Gamma_{j+2},\dotsc,\Gamma_k)\,.
\]
The domain of $\Gamma_j$ is $\Delta$. 
Then, for each $j=1,\dotsc,k$,  define the curve $C_j$ by  
\[
  C_j\ :=\ \{y\in\Delta\mid \psi_1(y)=\dotsb=\psi_{j-1}(y)
              \,=\,\Gamma_{j+1}(y)=\dotsb=\Gamma_k(y)\,=\,0\}\,.
\]
(The $j$th function $\psi_j$ or $\Gamma_j$ is omitted.)
This is a smooth curve by our assumption that the exponents $a_{i,j}$ are
sufficiently general.
Since we also assumed that the exponents are rational, we shall see that $C_j$
is even an algebraic curve.

If we iterate the Khovanskii-Rolle Theorem, we obtain
 \begin{eqnarray}
   |V(\psi_1,\psi_2,\dotsc,\psi_k)| &\leq& \flat(C_k)\,+\,   \nonumber
    |V(\psi_1,\psi_2,\dotsc,\psi_{k-1},\,\Gamma_k)|\\
  &\leq& \flat(C_k)\,+\,\flat(C_{k-1})\, +\,         \label{Eq:Kp_form}
     |V(\psi_1,\dotsc,\psi_{k-2},\,\Gamma_{k-1},\Gamma_k)|\\
  &\leq& \flat(C_k)\,+\,\dotsb\,+\,\flat(C_1)\, +\,
    |V(\Gamma_1,\Gamma_2,\dotsc,\Gamma_k)|\,.  \nonumber
 \end{eqnarray}
%

\begin{lemma}\label{L:Bounds}
 \ 
 \begin{enumerate}
  \item $\Gamma_{k-j}(y)\cdot\bigl(\prod_{i=1}^{n{+}k} p_i(y)\bigr)^{2^{j}}$ is a
     polynomial of degree $2^{j}\, n$.

  \item $C_j$ is a smooth algebraic curve and $\flat(C_j)\leq
  \frac{1}{2}\,2^{\binom{k-j}{2}}n^{k-j}\binom{n{+}k{+}1}{j}$.

  \item  $|V(\Gamma_1,\Gamma_2,\dotsc,\Gamma_k)|\leq 2^{\binom{k}{2}}n^k$.

 \end{enumerate}
\end{lemma} 

\begin{lemma}\label{L:technical}
 Let $k$ and $n$ be integers with $2\leq k,n$.
 Then
 \[
    \sum_{j=1}^k 2^{\binom{k-j}{2}}n^{k-j}\tbinom{n{+}k{+}1}{j}
    \ \leq \ \tfrac{e^2-1}{2} \cdot \, 2^{\binom{k}{2}}\, n^k\, .
 \]
\end{lemma}


\begin{proof}[Proof of Theorem~$\ref{Bounds}$]
 By~\eqref{Eq:Kp_form}, Lemma~\ref{L:Bounds} and Lemma~\ref{L:technical},
 $|V(\psi_1,\dotsc,\psi_k)|$ is at most
 \[
    2^{\binom{k}{2}} \, n^{k} + \frac{1}{2}
      \sum_{j=1}^k2^{\binom{k-j}{2}}n^{k-j} \tbinom{n+k{+}1}{j}
    \ <\ (1+\tfrac{e^2-1}{4})\, 2^{\binom{k}{2}} \, n^{k} \,,
 \]
 which proves Theorem~\ref{Bounds}.
\end{proof}

\begin{proof}[Proof of Lemma~{$\ref{L:technical}$}]
 Set $a_j := 2^{\binom{k-j}{2}}n^{k-j}\tbinom{n{+}k{+}1}{j}$ for
 $j=0,1,\dotsc,k$, 
 Note that $a_0=2^{\binom{k}{2}} n^k$.
 If $k=n=2$, then 
\[
  a_1+a_2\ =\ 20\ \leq\ \tfrac{e^2-1}{2} \cdot 8 \ =\ \tfrac{e^2-1}{2} \cdot a_0\,.
\]
(We have $3.1 < \tfrac{e^2-1}{2}< 3.2$.)
For other values of $(k,n)$, we will show that 
 \begin{equation}\label{Eq:compare}
   a_j\ \leq\ \frac{2^{j-1}}{j!} a_0\,.
 \end{equation}
 The lemma follows as 
 \begin{eqnarray*}
  \sum_{j=1}^k a_j &\leq& \
     \Bigl(\sum_{j=1}^k \frac{2^{j-1}}{j!} \Bigr)
      \cdot 2^{\binom{k}{2}}n^k \\
    &<& \Bigl(\sum_{j=1}^\infty \frac{2^{j-1}}{j!}\Bigr)
       \cdot 2^{\binom{k}{2}}n^k
    \ =\ \tfrac{e^2-1}{2}\cdot 2^{\binom{k}{2}}n^k\ .
 \end{eqnarray*}

 To show~\eqref{Eq:compare}, consider the expression
\[
   E(j,k,n)\ :=\ \frac{a_{j-1}}{a_j}\ =\ \frac{j \,n\, 2^{k-j}}{n{+}k{-}j{+}2}\,.
\]
 If we fix $j$ and $n$, then this is an increasing function of $k$ as 
 $(n{+}k{-}j{+}2)\cdot \ln 2>1$.
 Similarly, fixing $j$ and $k$, we obtain an increasing function of $n$,
 as $k{-}j{+}2>0$.

 Note that $E(j,j,2) = \frac{j}{2}$ and also that
\[
  E(1,2,3)\ =\ 1\quad\mbox{and}\quad
  E(1,3,2)\ =\ \tfrac{4}{3}\,.
\]
 Thus if $(k,n)\neq(2,2)$, we have $1\leq E(1,k,n)$ and so $a_1\leq a_0$ and
 if $j>1$, then $\frac{j}{2}=E(j,j,2) \leq E(j,k,n)$  as $j\leq k$ and $2\leq n$.
 Thus $a_j\leq \frac{2}{j} a_{j-1}$ for $j>1$ and $a_1 \leq a_0$. This leads to
 \[
   a_j\ \leq\ \frac{2^{j-1}}{j!} a_{1}\ \leq\ \frac{2^{j-1}}{j!} a_{0}
 \]
 which is~\eqref{Eq:compare}.
\end{proof}

We prove Lemma~\ref{L:Bounds}(1) and (2) in the following subsections.
Note that Lemma~\ref{L:Bounds}(3) follows from Lemma~\ref{L:Bounds}(1).
Since $\prod_{i=1}^{n{+}k} p_i(y)$ does not vanish on $\Delta$, if we set
 \begin{equation}\label{Eq:Poly_form_of_Gamma}
   F_{j}\ :=\ \Gamma_{j}\cdot \Bigl(\prod_{i=1}^{n{+}k} p_i(y)\Bigr)^{2^{k-j}}\,,
 \end{equation}
 then by Lemma~\ref{L:Bounds}(1) this is a polynomial of degree $2^{k-j} n$ and
\[
   V(F_1,F_2,\dotsc,F_k)\ =\ V(\Gamma_1,\Gamma_2,\dotsc,\Gamma_k)\,.
\]
 But then B\'ezout's Theorem gives the upper bound
\[
  \deg(F_1)\deg(F_2)\dotsb\deg(F_k)\ =\ 
   2^{(k{-}1) + \dotsb +2+1+0}n^k \ =\ 2^{\binom{k}{2}}n^k
\]
 for the number of isolated common zeroes of $F_1,\dotsc,F_k$ in $\R^k$, which
 is a bound for $V(\Gamma_1,\dotsc,\Gamma_k)$, the common zeroes of
 $\Gamma_1,\dotsc,\Gamma_k$ in $\Delta$.

\subsection{Proof of Lemma~\ref{L:Bounds}(1)}
 We first give a useful lemma, and then
 determine the form of the functions $\Gamma_i(y)$.

 Let $m\geq k$ be integers.
 Then $\binom{[m]}{k}$ is the collection of subsets of $\{1,\dotsc,m\}$
 with cardinality $k$.
 Given a vector $c=(c_1,\dotsc,c_m)$, and $I=\{i_1<\dotsb<i_k\}\in\binom{[m]}{k}$,
 set $c_I:=c_{i_1}\dotsb c_{i_k}$.
 If $D$ is a $m\times k$ matrix and $I\in\binom{[m]}{k}$, 
 let $D_I$ be the determinant of the $k\times k$ submatrix of $D$ formed by the
 rows indexed by $I$. 

\begin{lemma}\label{L:useful}
 Let $c=(c_1,\dotsc,c_m)$ be a vector of indeterminates
 and $D=(d_{i,j})$ and $E=(e_{i,j})$ be $m\times k$ matrices of
 indeterminates.
 Then 
\[
   \det\Bigl( \sum_{i=1}^m  c_i d_{i,j} e_{i,l} \Bigr)_{j,l=1,\dotsc,k}
   \ =\ \sum_{I\in\binom{[m]}{k}}c_I  D_I E_I \,.
\]
\end{lemma}
\begin{proof}
This follows from the Cauchy-Binet formula.
\end{proof}

We use Lemma~\ref{L:useful} to determine the form of the functions $\Gamma_j$.
Since 
\[
  \frac{\partial\psi_j}{\partial y_l}\ =\ 
   \sum_{i=1}^{n+k}\frac{a_{i,j} b_{i,l}}{p_i}\,,
\]
Lemma~\ref{L:useful} implies that 
 \begin{equation}\label{E:Gamma_k}
   \Gamma_k\ =\ \sum_{I\in\binom{[n+k]}{k}} \frac{A_I B_I}{p_I}\,.
 \end{equation}

Lemma~\ref{L:Bounds}(1) is a consequence of the following computation.
The set $\binom{[n+k]}{k}^{2^j}$ consists of lists $\calI$ of $2^j$ elements of
$\binom{[n+k]}{k}$. 
For each such list $\calI \in \binom{[n+k]}{k}^{2^j}$, set
\[
   B(\calI)\ :=\ \prod_{I\in \calI} B_I 
   \qquad\mbox{and}\qquad
   p(\calI)\ :=\ \prod_{I\in \calI} p_I\,.
\] 

\begin{lemma}\label{L:Gen_Gamma}
 There exist numbers $A(\calI)$ for each $\calI\in\binom{[n+k]}{k}^{2^j}$, 
 so that 
 \begin{equation}\label{Eq:Gen_Gamma}
   \Gamma_{k-j}\ =\ \sum_{\calI\in\binom{[n+k]}{k}^{2^j}}
     \frac{A(\calI) B(\calI)}{p(\calI)}\,.
 \end{equation}
\end{lemma}

\begin{proof}
%
%
 We prove this by induction on $j$.
 The case  $j=0$ is formula~\eqref{E:Gamma_k}.
 We have
\[
  \Gamma_{k-j}\ :=\ J(\psi_1,\dotsc,\psi_{k-j},
     \Gamma_{k-j+1},\dotsc,\Gamma_k)\,.
\]
 Since the Jacobian is multilinear, we use the formulas~\eqref{Eq:Gen_Gamma}
 for $\Gamma_{k-j+1},\dotsc,\Gamma_k$ to obtain
\[
  \Gamma_{k-j}\ =\ 
   \sum_{\calI^{(j-1)}\in\binom{[n+k]}{k}^{2^{j-1}}}\ \dotsb\ 
         \sum_{\calI^{(0)}\in\binom{[n+k]}{k}^{2^0}}\ 
        \prod_{r=0}^{j-1} A(\calI^{(r)})B(\calI^{(r)})
    \ \cdot\   J(\calI^{(j-1)},\dotsc,\calI^{(0)})\,,
\]
where $J(\calI^{(j-1)},\dotsc,\calI^{(0)}):=
J\bigl(\psi_1,\dotsc,\psi_{k-j},
      \frac{1}{p(\calI^{(j-1)})},\,\dotsc,\,
       \frac{1}{p(\calI^{(0)})}\bigr)$, which equals
 \[     \biggl(\prod_{r=0}^{j-1}\frac{1}{p(\calI^{(r)})}\biggr) \cdot
   J\Bigl(\psi_1,\dotsc,\psi_{k-j},\,
            \log\Bigl(\tfrac{1}{p(\calI^{(j-1)})}\Bigr),
    \dotsc, \log\Bigl(\tfrac{1}{p(\calI^{(0)})}\Bigr)\Bigr)\,.
 \]
%

Since
\[
   \frac{\partial}{\partial y_l}\log\Bigl(\frac{1}{p(\calI^{(r)})}\Bigr)
    \ =\ 
   \sum_{i=1}^{n+k} \frac{a'_{i,r} b_{i,l}}{p_i(y)}\,,
\]
where $-a'_{i,l}$ counts the number of sets 
$I$ in the list $\calI^{(l)}$ which contain $i$, Lemma~\ref{L:useful}
implies that $J(\calI^{(j-1)},\dotsc,\calI^{(0)})$ has the form
 \begin{equation}\label{Eq:new_form}
   \sum_{I\in\binom{[n+k]}{k}} \frac{A'_I B_I}{p_I}\,,
 \end{equation}
where $A'$ is the matrix whose entry in column $r$ is
$a_{i,r}$ if $r\leq k-j$ and is $a'_{i,k-r}$ otherwise.
If we substitute~\eqref{Eq:new_form} into the previous formulas and write
it as a single sum, we obtain an expression of the form~\eqref{Eq:Gen_Gamma}.
\end{proof}

%
\subsection{Proof of Lemma~\ref{L:Bounds}(2)}
%
 Let $\mu_j\subset\Delta$ be defined by the equations
 \[
   \psi_1(y)\ =\ \psi_2(y)\ =\ \dotsb\ =\ \psi_{j-1}(y)\ =\ 1\,.
 \]

 The polyhedron $\Delta$ is an open subset of $\R^k$ and we let
 $\overline{\Delta}$ be its closure in the real projective space $\R\P^k$.
 If $\Delta$ is bounded, then $\overline{\Delta}$ is its closure, which is a
 polytope. 
 If $\Delta$ is unbounded, then $\overline{\Delta}$ is combinatorially
 equivalent to the polytope
 \begin{equation}\label{Eq:poly}
   \Delta\cap\{y\mid v\cdot y\leq r\}\,,
 \end{equation}
if $v+\Delta\subset \Delta$ and $r$ is large enough so that all the vertices of
$\Delta$ lie in~\eqref{Eq:poly}.

 A \Blue{{\it PL-manifold}} is a manifold in the category of piecewise-linear
 spaces. 
 Given a polyhedral complex $P$, let $M_j(P)$ be the maximum number of 
 $j$-dimensional faces ($j$-faces) whose union forms a $j$-dimensional
 PL-submanifold of $P$. 

\begin{lemma}\label{L:mu}
 The set $\mu_j$ is a smooth algebraic subset of $\Delta$ of
 dimension $k{+}1{-}j$. 
 The intersection of its closure $\overline{\mu}_j$ with the boundary of
 $\overline{\Delta}$ is a union of\/
 $(k{-}j)$-faces of\/ $\partial\overline{\Delta}$ which forms a
 PL-submanifold of $\partial\overline{\Delta}$.
\end{lemma}

 Since $C_j$ is the subset of $\mu_j$ cut out by the 
 polynomials $F_{j+1},\dotsc,F_k$~\eqref{Eq:Poly_form_of_Gamma}, it is an
 algebraic subset of $\Delta$ and, as we already noted, a smooth curve.
 Our choice of generic exponents $a_{i,l}$ ensures that each  unbounded
 component of $C_j$ meets the boundary of $\Delta$ in two distinct points, and
 no two components meet the boundary at the same point.

\begin{cor}\label{C:boundary_bounds}
\qquad$\flat(C_j)\ \leq\ 
      \frac{1}{2}\,2^{\binom{k-j}{2} n^{k-j}}\,
      M_{k-j}(\overline{\Delta})$.
\end{cor}

 Corollary~\ref{C:boundary_bounds}  implies the bound of
 Lemma~\ref{L:Bounds}(2).
 Indeed, if $\Phi_j(\overline{\Delta})$ is the number of $j$-faces of
 $\overline{\Delta}$, then 
 \begin{equation}\label{Eq:Compare}
   M_{k-j}(\overline{\Delta})\ \leq\ \Phi_{k-j}(\overline{\Delta})
   \ \leq\ \tbinom{n+k+1}{j}\,.
 \end{equation}
 The last inequality follows as each $(k{-}j)$-face lies on $j$ of the
 $n{+}k{+}1$ facets.

\begin{proof}[Proof of Corollary~$\ref{C:boundary_bounds}$]
 Let $\overline{C_j}$ be the closure of $C_j$ in $\overline{\Delta}$.
 First, 
 $\flat(C_j) =\frac{1}{2}|\overline{C}_j\cap \partial\overline{\Delta}|$
 as each unbounded component of $C_j$ contributes two points to
 $\overline{C}_j\cap\partial\overline{\Delta}$.
 The points in $\overline{C}_j\cap\partial\overline{\Delta}$ are
 points of $\overline{\mu}_j\cap\partial\overline{\Delta}$ where the
 polynomials $\overline{F}_{j+1},\dotsc,\overline{F}_k$ vanish.
 For each $(k{-}j)$-dimensional face $\phi$ in
 $\overline{\mu}_j\cap\partial\overline{\Delta}$, 
 there will be at most 
\[
   2^{\binom{k-j}{2}}n^{k-j}\ =\ 
    \deg(\overline{F}_{j+1})\cdots \deg(\overline{F}_{k-1})
     \cdot\deg(\overline{F}_{k})
\]
 points of $\overline{C}_j\cap\phi$.
 Indeed, this is the B\'ezout bound for the number of isolated common zeroes of the
 polynomials $\overline{F}_{j+1},\dotsc,\overline{F}_k$ on the affine span of
 the face $\phi$.
 Since these faces $\phi$ form a PL-submanifold of $\partial\overline{\Delta}$,
 we deduce the bound of the corollary.
\end{proof}

\begin{proof}[Proof of Lemma~$\ref{L:mu}$]
 Since $\psi_i$ and $f_i$ define the same subset of $\Delta$, $\mu_j$ is defined
 on $\Delta$ by 
 \begin{equation}\label{Eq:mu_def}
   f_1(y)\ =\ f_2(y)\ =\ \dotsb\ =\ f_{j-1}(y)\ =\ 1\,.
 \end{equation}
 The exponents $a_{i,l}$ are rational numbers.
 Let $N$ be their common denominator.
 As the functions $f_i$ are positive on $\Delta$, we see that
 $\mu_j\subset\Delta$ is defined by the rational functions
\[
  f_1^N(y)\ =\ f_2^N(y)\ =\ \dotsb\ =\ f_{j-1}^N(y)\ =\ 1\,,
\]
 and is thus algebraic.
 Since the exponents are also general,
 we conclude that $\mu_j$ is smooth and has the expected dimension $k{-}j{+}1$.

 Our arguments to prove the remaining statements in Lemma~\ref{L:mu} are
 local, and in particular, we will  work 
 near points on the
 boundary of $\overline{\Delta}\subset\R\P^k$.
 For this, we may need to homogenize the functions $f_i$ and $F_i$.
 We first homogenize the linear polynomials $p_i(y)$ with respect to a
 new variable $y_0$,
\[
   \overline{p}_i(y)\ :=\ 
    b_{i,0}y_0 + b_{i,1}y_1 + \dotsb + b_{i,k}y_k\,.
\]
 Let $\overline{F}_j$ be the homogeneous version of $F_j$, where we replace
 $p_i$ by $\overline{p}_i$.
 For the functions $f_j$, we first replace $p_i$ by $\overline{p}_i$ and then
 multiply by an appropriate power of $y_0$,
\[
  \overline{f}_j\ := 
   y_0^{-\sum_i a_{i,j}}\prod_{i=1}^{n+k} \overline{p}_i(y)^{a_{i,j}}\,.
\]

We show that the intersection of $\mu_j$ with the boundary of
$\overline{\Delta}$ is a union of faces of $\overline{\Delta}$ of dimension
$k{-}j$.
For this, select a face $\phi$ of $\overline{\Delta}$ of dimension
$k{-}j$.
Changing coordinates in $\R\P^k$ if necessary, we may assume that 
the face lies along the coordinate plane
defined by $y_1=\dotsb=y_j=0$.
We also reindex the forms $p_i$ so that $p_i=y_i$ for $i=1,\dotsc,k$.

Since every square submatrix of the matrix $A$ is
invertible, we may suppose that its first $j{-}1$ rows and columns
form the indentity matrix---without changing the column span of
its first $j{-}1$ columns.
By Remark~\ref{R:lin=mult}, 
if we multiply every entry in $A$ by $-1$, we may assume that 
the equations defining $\mu_j$ have the form
 \begin{equation}\label{Eq:New_mu}
   y_i\ =\  y_{j}^{a_{j,i}}\cdot\prod_{l=j+1}^{n+k} p_l(y)^{a_{l,i}}\,,
   \qquad\mbox{for}\quad i=1,\dotsc,j{-}1\,,
 \end{equation}
and no exponent $a_{j,i}$ vanishes.

We first observe that the face $\phi$ lies in the closure $\overline{\mu}_j$ of
$\mu_j$ if and only if the exponents $a_{j,i}$ for $i=1,\dotsc,j{-}1$ are all
positive.
Moreover, $\overline{\mu}_j$ cannot contain any point in the relative interior
of a face which properly contains $\phi$, as on that face we still have
$p_l(y)>0$ for $l>j$, at least one of the variables $y_1,\dotsc,y_j$ will 
vanish, and at least one is non-vanishing.
But if one of the variables $y_1,\dotsc,y_j$ vanishes, then the
equations~\eqref{Eq:New_mu} imply that they all do, which is a contradiction.

Now consider the variety $\mu_j$ in the neighborhood of a vertex of the face
$\phi$. 
Applying a further affine linear change to our variables $y$, we may assume this
vertex lies at the origin of $\R^k$.
Then the equations defining $\mu_j$ have the form
\[
   y_i\ =\  y_j^{a_{j,i}} y_{j+1}^{a_{j+1,i}}\dotsb y_k^{a_{k,i}}\cdot
     \prod_{l=k+1}^{n+k} p_l(y)^{a_{l,i}}\,,
   \qquad\mbox{for}\quad i=1,\dotsc,j{-}1\,.
\]
Since $p_l(0)>0$ for $l>k$, we see that $\mu_j$ is approximated in this
neighborhood of the origin by the zero set of equations of the form
\[
   y_i\ =\  y_j^{a_{j,i}} y_{j+1}^{a_{j+1,i}}\dotsb y_k^{a_{k,i}}\cdot b_i\,,
   \qquad\mbox{for}\quad i=1,\dotsc,j{-}1\,.
\]
where $b_i>0$.
We may scale $y_i$ so that $b_i=1$, and then these equations define a 
subtorus of the positive part of the real torus $\R^k_>$.
Since the exponents are rational numbers, the closure of this subtorus in the
non-negative orthant is the non-negative part of a (not necessarily normal)
toric variety, which is homeomorphic to a polytope~\cite[\S~4]{Fu93}.
In particular, its boundary is homeomorphic to the boundary of a polytope and is
thus a manifold.

This proves that, in a neighborhood of any vertex, the intersection of
$\overline{\mu_j}$ with the boundary of $\Delta$ forms a manifold.
But this implies that the union of faces of $\overline{\Delta}$ which are
contained in $\overline{\mu}_j$ forms a PL-manifold, and completes the proof of
the lemma.
\end{proof}

\subsection{Concrete bounds for $k=2, 3$}
When $k=2$, the estimate~\eqref{Eq:Kp_form} becomes
\[
   |V(\psi_1, \psi_2)| \ \leq\ \flat(C_1)\ +\ \flat(C_2)\ +\ 
   |V(\Gamma_1,\Gamma_2)|\,.
\]

\begin{thm}\label{T:Boundssolfork=2}
 Suppose that $k=2$.
 Then
\begin{enumerate}
 \item $\flat(C_2) \leq \lfloor\frac{n+3}{2}\rfloor$, 

 \item $\flat(C_1) \leq \lfloor\frac{n(n+3)}{2}\rfloor$, and \rule{0pt}{13pt}

 \item $|V(\Gamma_1,\Gamma_2)|\leq 2n^2$.\rule{0pt}{13pt}
\end{enumerate}
Thus a fewnomial system of $n$ polynomials in $n$
variables having $n+3$ monomials has at most
$2n^2+\lfloor\frac{(n+1)(n+3)}{2}\rfloor$ positive solutions. 
\end{thm}

\begin{proof}
 $\overline{\Delta}$ is a polygon with at most
 $n+3$ edges and $n+3$ vertices. 
 By Corollary~\ref{C:boundary_bounds}  and~\eqref{Eq:Compare}, we have
 $2\flat(C_2) \leq \Phi_0(\overline{\Delta}) \leq n+3$
 and $2\flat(C_1) \leq n\Phi_1(\overline{\Delta}) \leq n(n+3)$.
 By Lemma~\ref{L:Bounds}(3), we have $|V(\Gamma_2,\Gamma_1)|\leq 2n^2$.
\end{proof}

In particular, when $n=k=2$, this gives a bound of $15$ for the number of
positive solutions to a fewnomial system of $2$ polynomials in $2$ variables and
with $5$ monomials. \smallskip 

When $k=3$, the estimate~\eqref{Eq:Kp_form} becomes
%
\[
   |V(\psi_1, \psi_2,\psi_3)| \ \leq\ \flat(C_1)\ 
    +\ \flat(C_2)\ +\ \flat(C_3)\ +\
   |V(\Gamma_1,\Gamma_2, \Gamma_3)|\,.
\]

\begin{thm}\label{T:Boundssolfork=3}
 Suppose that $k=3$.
 Then
\begin{enumerate}
\item $\flat(C_3) \leq n+2$, 

\item $\flat(C_2) \leq n(n+2)$,\rule{0pt}{13pt}
 
\item $\flat(C_1) \leq n^2(n+4)$,\rule{0pt}{13pt}

\item $|V(\Gamma_1,\Gamma_2,  \Gamma_3)|\leq 8n^3$.\rule{0pt}{13pt}
\end{enumerate}
Thus 
a fewnomial system of $n$ polynomials in $n$
variables having 
$n+4$ 
monomials has at most $9n^3+5n^2+3n+2$ positive solutions.
\end{thm}

\begin{proof}
 The polytope $\overline{\Delta}$ is a simple three-dimensional polytope
 and Corollary~\ref{C:boundary_bounds} gives

\begin{enumerate}
\item $2\flat(C_3) \leq M_0(\overline{\Delta})$,

\item $2\flat(C_2) \leq n M_1(\overline{\Delta})$,\rule{0pt}{13pt}

\item $2\flat(C_1) \leq 2n^2M_2(\overline{\Delta})$.\rule{0pt}{13pt}
\end{enumerate}

 We have $M_j(\overline{\Delta}) \leq \Phi_j(\overline{\Delta})$
 and McMullen's Upper Bound Theorem~\cite{Mc70} applied to the dual of our
 polytope gives 
 $\Phi_0(\overline{\Delta}) \leq 2(n+2)$, $\Phi_1(\overline{\Delta}) \leq 3(n+2)$, 
 and $\Phi_2(\overline{\Delta}) \leq n+4$.
 This gives estimates for $M_j(\overline{\Delta})$. 
 However, $M_1(\overline{\Delta})\leq \Phi_0(\overline{\Delta})$, as a
 1-dimensional PL-manifold is a union of closed polygons. 
 Combining these bounds proves the result.
\end{proof}

%
\section{Number of compact components of hypersurfaces}
%

Let $f$ be a polynomial with $n{+}k{+}1$ monomials whose real exponents 
affinely span $\R^n$. 
Let $V(f)\subset \R^n_>$ be the set of zeroes of $f$ in the
positive orthant, which we assume is smooth, and let $\kappa(f)$ be the number
of compact components of $V(f)$. 
When $n=1$, a bound of $n{+}k$ for $\kappa(f)$ is given by Descartes's rule of
signs~\cite{D1636}, and when $k=1$ the bound of $1$ for $\kappa(f)$ is proven by
Bihan, Rojas, and Stella~\cite{BRS}.
We therefore assume that $n,k \geq 2$.

\begin{thm}\label{Th:Components}
 \ 
 \begin{enumerate}
  \item We have $\kappa(f)<\frac{e^2+3}{8}\cdot 2^{\binom{k}{2}} n^k$.
  \item We also have $\kappa(f)< C(k) n^{k-1}$,
        where $C(k)$ has order $O(k2^{\binom{k}{2}})$.
   \end{enumerate}
\end{thm}

The first assertion is a straightforward application of
Theorem~\ref{FinalBounds}.
When $k>n$, it implies the second assertion.
The proof of the second assertion when $k\leq n$ exploits the sparsity of a 
Gale system.

After an analytic change of variables, we may assume that $f$ has the
form
\[
  f\ =\ \sum_{i=1}^n e_i z_i \ +\ c_1 z^{a_1} \ + \cdots +\ c_k z^{a_k}\ +\ e_0\,,
\]
where $e_i\in\{\pm1\}$, $a_j\in\R^n$ and the coefficients $c_i$ are non-zero
Its support is $\calA:=\{0,a_1,\dotsc,a_k,\, e_1,\dotsc,e_n\}$.

Since $V(f)$ is smooth, the coordinate function $z_1$ has at least
2 critical points on each bounded component of $V(f)$.
These critical points are solutions to the system
 \begin{equation}\label{Eq:Crit_Point}
   f\ =\ z_2\frac{\partial f}{\partial z_2}\ =\ 
         z_3\frac{\partial f}{\partial z_3}\ =\ 
          \dotsb \ =\ 
         z_n\frac{\partial f}{\partial z_n}\ =\ 0\,.
 \end{equation}
(Since we work in $\R^n_>$, no coordinate $z_i$ vanishes.)
The resulting system has the form
 \begin{equation}\label{Eq:comp_sys}
   \begin{array}{rcl}
        -e_1z_1&=& 
       {\displaystyle\sum_{i=2}^n} e_i z_i \ +\ c_1 z^{a_1} \ 
      + \cdots + \ c_k z^{a_k}\ +\ e_0\\
     - e_jz_j&=& c_1a_{1,j}z^{a_1} + \cdots + \ c_k a_{k,j}z^{a_k}
    \qquad\qquad\mbox{for}\quad j=2,3,\dotsc,n\,.
   \end{array}
 \end{equation}
This is a system of $n$ polynomials in $n$ variables with support $\calA$,
which we may assume has only non-degenerate solutions by perturbing the
parameters without altering the topology of $V(f)$, as $V(f)$ is smooth.
Theorem~\ref{Th:Components}(1) follows now from Theorem~\ref{FinalBounds}
as $2\kappa(f)$ is bounded by the number of solutions to~\eqref{Eq:comp_sys}.

By Theorem~\ref{Th:GaleSystem}, the number of solutions in $\R^n_>$
to~\eqref{Eq:Crit_Point} is equal to the number of solutions to an associated
Gale system in a polyhedron $\Delta\subset\R^k$.

Performing Gaussian elimination on~\eqref{Eq:comp_sys}, we obtain a system of
the form 
 \begin{equation}\label{Eq:component_eqs}
  \begin{array}{rclcl}
   z_1&=& b_{1,1}z^{a_1} + \dotsb + b_{1,k}z^{a_k} +b_{1,0}

      &=:&p_1(z^{a_1}, \dotsc, z^{a_k})\\
   z_2&=& b_{2,1}z^{a_1} + \dotsb + b_{2,k}z^{a_k}

      &=:&p_2(z^{a_1}, \dotsc, z^{a_k})\\
    & \vdots&&\vdots&\\

 z_n&=& b_{n,1}z^{a_1} + \dotsb + b_{n,k}z^{a_k}

      &=:&p_n(z^{a_1}, \dotsc, z^{a_k})\\
  \end{array}
 \end{equation}
where $b_{i,j} \in \R$.

Set $y_i:=z^{a_i}$, for $i=1,\dotsc,k$, so that the right
hand sides of~\eqref{Eq:component_eqs} are linear functions $p_i(y)$ of
$y$. 
Then use the exponents $a_1, \dotsc, a_k$ to get the Gale system
\begin{equation}\label{Eq:Comp_Gale}
   \Blue{f_j}(y)\ :=\ y_j^{-1}\cdot(p_1(y))^{a_{1,j}} \cdot 
    \prod_{i=2}^n p_i(y)^{a_{i,j}} \ =\ 1\qquad 
     \mbox{for}\quad j=1,2,\dotsc,k \,.
 \end{equation}
This has the form of~\eqref{Eq:GaleSystem}, but the linear factors
$p_i$ for $i>1$ are sparse. 
As in Section~\ref{S:general_k}, we may assume that
every square submatrix of $A=(a_{i,j})$ is 
invertible, and that the exponents $a_{i,j}$ are rational numbers.

\begin{thm}\label{T:boundscomponent}
 Let $f$ be a polynomial with real exponents in $n$ variables having $n{+}k{+}1$  
 monomials with $k\leq n$.
 Then
 \[
   \kappa(f)\ \leq\ 
   \left(\tfrac{k}{2}2^{\binom{k}{2}}\,+\,
     \tfrac{e^2+1}{8} \cdot k 2^{\binom{k-1}{2}} \right) n^{k-1}
      \ +\ 
    \left( \tfrac{e^2}{8} \cdot 2^{\binom{k-2}{2}} \right) n^{k-2}\,.
 \]
\end{thm}

This implies Theorem 4.1(2).

We prove Theorem~\ref{T:boundscomponent} with the methods of
Section~\ref{S:general_k}, taking into account the 
sparsity of the Gale system~\eqref{Eq:Comp_Gale}. 
Consider the polyhedron $\Delta\subset\R^k_>$ defined by
 \begin{equation}\label{Eq:Delta_def}
  \Delta\ :=\ \{ y\in\R^k_>\mid p_i(y)>0\quad\mbox{for}\quad
      i=1,\dotsc,n\}\,.
 \end{equation}
For each $j=1,\dotsc,k$,  set $p_{n+j}=y_j$.
Since the linear polynomials $p_2,\dotsc,p_{n+k}$ do not have a constant term
but are otherwise general, they define a cone with vertex the origin and base 
a simple polytope $\phi^{\infty}$ with at most $n{+}k{-}1$ facets lying in the
hyperplane at infinity.
We obtain $\Delta$ from this cone by adding the single inequality $p_1(y)>0$.
This gives at most one new facet, $\phi^1$.
By our assumptions on the genericity of the polynomials $p_i$, every face of
dimension $k{-}j$ lies on $j$ facets, unless $j=k$ and that face is the origin.

Let $\overline{\Delta}\subset\R\P^k$ be the closure of $\Delta$.
Let $\Phi_j^\ell(\overline{\Delta})$ count the \Blue{{\it linear faces}} of
dimension $j$, those whose affine span contains the origin,
and let $\Phi_j^{n\ell}(\overline{\Delta})$ count the $j$-faces of 
$\overline{\Delta}$ lying on either $\phi^\infty$ or $\phi^1$.

For each $j=1,\dotsc,k$, let $\psi_j$ be the logarithm of right hand side
of the corresponding equation in~\eqref{Eq:Comp_Gale}.
Define functions $\Gamma_k(y),\dotsc,\Gamma_1(y)$ by the recursion 
\[
   \Gamma_j\ :=\ J(\psi_1,\psi_2,\dotsc,\psi_j,\, 
              \Gamma_{j+1},\Gamma_{j+2},\dotsc,\Gamma_k)\,.
\]
Then, for $j=1,\dotsc,k$,  define the curve $C_j$ by  
\[
  C_j\ :=\ \{y\in\Delta\mid \psi_1(y)=\dotsb=\psi_{j-1}(y)
              \,=\,\Gamma_{j+1}(y)=\dotsb=\Gamma_k(y)\,=\,0\}\,.
\]
Successive uses of the Khovanskii-Rolle theorem give an upper bound
for $\kappa(f)$,
 \begin{equation}\label{E:critical}
   \kappa(f)\ \leq\ 
    \tfrac{1}{2}\bigl(|V(\psi_1,\dotsc,\psi_k)|\ \leq\ 
     \tfrac{1}{2}\bigl(\flat(C_k)\,+\,\dotsb\,+\,\flat(C_1)\, +\,
    |V(\Gamma_1,\Gamma_2,\dotsc,\Gamma_k)|\bigr)\,.
 \end{equation}

For $j=1,\dotsc,k$, define $F_j$ by formula~\eqref{Eq:Poly_form_of_Gamma}.
This polynomial of degree $n2^{k-j}$ defines the same subset of
$\Delta$ as does $\Gamma_j$.
The sparsity of $p_i$ (for $i>2$) implies the sparsity of 
polynomials $F_j$.

\begin{lemma}\label{L:Boundscomp}
 \ 
 \begin{enumerate}

\item Every monomial of $F_{j}$ has degree in the interval 
$[(n{-}1)2^{k-j}  ,  n2^{k-j}] $.
%

\item $C_j$ is a smooth algebraic curve, 
    $\flat(C_k)\leq \frac{1}{2}\bigl(1+\Phi_{0}^{n\ell}(\overline{\Delta})\bigr)$,  
    and for $j=1,\dotsc,k-1$ we have 
\[
    \flat(C_j) \  \leq \   \tfrac{1}{2}\Bigl(
         2^{\binom{k-j}{2}} \, n^{k-j} 
               \Phi_{k-j}^{n\ell}(\overline{\Delta})
              \  +\   
     2^{\binom{k-j}{2}} \,\left( n^{k-j} -{(n{-}1)}^{k-j} \right) 
         \, \Phi_{k-j}^\ell(\overline{\Delta}) \Bigr)\,.  
\]
\item  $|V(\Gamma_1,\Gamma_2,\dotsc,\Gamma_k)| \leq 
              2^{\binom{k}{2}} \, \left( n^k- {(n{-}1)}^k \right)$.
 \end{enumerate}
\end{lemma}

Lemma~\ref{L:Boundscomp} is proved in the next sections.

\begin{lemma}\label{L:facescompact}
We have
\begin{enumerate}

\item $\Phi_{k-j}^{n\ell}(\overline{\Delta}) \leq 
                    2\binom{n+k-1}{j-1}+\binom{n+k-1}{j-2}$,
       for $j=1,\dotsc,k$,

\item $\Phi_{k-j}^\ell(\overline{\Delta}) \leq \binom{n+k-1}{j}$,
       for $j=1,\dotsc,k{-}1$, \rule{0pt}{14pt}
       and $\;\Phi^\ell_0(\overline{\Delta})=1$.

\end{enumerate}
\end{lemma}

\begin{proof}
 The faces of $\overline{\Delta}$ arise in (possibly) three different ways.
 A $(k{-}j)$-face $\phi$ either $(a)$ lies in the intersection of
 the facets $\phi^\infty$ and $\phi^1$, or it $(b)$ lies in exactly one of
 $\phi^\infty$ or $\phi^1$, or else $(c)$ its affine span contains the origin.
 Unless $\phi$ is the origin, it lies on $j$ facets.
 Among those are some \Blue{{\it linear facets}} which are defined by one of
 $p_2,\dotsc,p_{n{+}k}$.
 In case $(a)$ there are $j{-}2$ such facets, in case $(b)$ there  are $j{-}1$
 facets, and in case $(c)$ all $j$ facets are linear.  
 We may bound the $(k{-}j)$-faces of each type by enumerating the
 possible sets of linear facets on which they lie,
\[
   (a)\ \  \binom{n{+}k{-}1}{j{-}2}\,,\qquad
   (b)\ \ 2\binom{n{+}k{-}1}{j-1}\,,
     \quad\mbox{and}\quad
   (c)\ \  \binom{n{+}k{-}1}{j}\,.
\]
 The factor of 2 in the estimate $(b)$ is because such faces lie on either
 $\phi^\infty$ or $\phi^1$.

 The lemma follows from this.
 For (1), note that $\Phi_{k-j}^{n\ell}(\overline{\Delta})$ 
 counts faces of types $(a)$ and $(b)$.
 For (2), $\Phi_{k-j}^\ell(\overline{\Delta})$ counts faces of type $(c)$, when
 $k{-}j>0$. 
\end{proof}

\begin{proof}[Proof of Theorem~$\ref{T:boundscomponent}$]
 Substitute the estimates of Lemma~\ref{L:Boundscomp} into~\eqref{E:critical} to
 obtain
 \begin{multline*}
    \kappa(f)\ \leq\ \frac{1}{4}\ +\ 
    \frac{1}{4}\sum_{j=1}^k 2^{\binom{k-j}{2}}n^{k-j} 
       \Phi^{n\ell}_{k-j}(\overline{\Delta})\\
    \ +\ \frac{1}{4}\sum_{j=1}^k 2^{\binom{k-j}{2}}
              (n^{k-j}-(n{-}1)^{k-j}) \Phi^\ell_{k-j}(\overline{\Delta})
        \ +\ \frac{1}{2} 2^{\binom{k}{2}}(n^k-(n{-}1)^k)\,.
 \end{multline*}
(We used $2^{\binom{k-k}{2}}n^{k-k}=1$ in the first sum and
 $n^{k-k}-(n{-}1)^{k-k}=0$ in the second.)
Substituting the estimates of Lemma~\ref{L:facescompact} and rearranging, we
 obtain the estimate
\[
   \kappa(f)\ \leq\ \frac{1}{4}\sum_{j=1}^k   A_j\ \,+\ \,
     \frac{1}{4}\sum_{j=2}^k B_j\ \, +\ \,  \frac{1}{4}\ +\ 
     \frac{1}{2}2^{\binom{k}{2}} (n^k-(n-1)^k)\,,
\]
where 
\[
   A_j\ :=\ 2^{\binom{k-j}{2}} 
    \left(\left(n^{k-j}-(n-1)^{k-j}\right) \tbinom{n+k-1}{j} 
    + 2n^{k-j} \tbinom{n+k-1}{j-1} \right)\, ,
 \]
 and 
 \[
 B_j\ :=\ 2^{\binom{k-j}{2}} n^{k-j}\tbinom{n+k-1}{j-2}\,.
 \]

 We have $n^{k-j}-(n-1)^{k-j} \leq (k-j)n^{k-j-1}$, as the function
 $x^{k-j}$ is (weakly) convex and non-decreasing.
 Thus 
 \[
   A_j\ \leq\  2^{\binom{k-j}{2}} n^{k-j}
         \left( \tfrac{k-j}{n}\tbinom{n+k-1}{j} +2\tbinom{n+k-1}{j-1} \right)\,.
 \]
 Since $1\leq j\leq k\leq n$, we have
 \[
  \tfrac{k{-}j}{n} \tbinom{n{+}k{-}1}{j}\ =\ 
    \tfrac{k-j}{j} \, \tfrac{n{+}k{-}j}{n}\,
        \tfrac{n{+}k{-}j{+}1}{n{+}k} \tbinom{n{+}k}{j{-}1}
    <\  2\tfrac{k-j}{j}\tbinom{n{+}k}{j{-}1}\ .
\]
But we also have $\binom{n{+}k{-}1}{j{-}1}\leq\binom{n{+}k}{j{-}1}$, and so
 \[
   A_j\ <  2\cdot 2^{\binom{k-j}{2}} n^{k-j}
         \left( \tfrac{k-j}{j}\tbinom{n{+}k}{j{-}1} +\tbinom{n{+}k}{j{-}1}\right)\ 
      =\  2k\cdot 2^{\binom{k-j}{2}} n^{k-j}\tbinom{n{+}k}{j{-}1}\,.
 \]
 Thus 
 \[
   \sum_{j=1}^k A_j\ \leq\ 
     2k\sum_{j=1}^k 2^{\binom{k-j}{2}}\tfrac{n^{k-j}}{j}\tbinom{n{+}k}{j{-}1}\ .
 \]
 If we set $i:=j-1$ and $l:=k{-}1$, then we obtain
 \[
  \sum_{j=1}^k A_j\ \leq\ 
    2k\Bigl( 2^{\binom{k-1}{2}}\, n^{k-1}\ +\ 
    \sum_{i=1}^{l+1} 2^{\binom{l-i}{2}}
    \tfrac{n^{l-i}}{i{+}1}\tbinom{n{+}l{+}1}{i}\Bigr)\ .
 \]
 Arguing as in Lemma~\ref{L:technical} bounds the second sum by
\[
   2^{\binom{l}{2}} n^l\,\sum_{i=1}^\infty \tfrac{2^{i-1}}{(i+1)!}
   \ =\ \tfrac{e^2-3}{4}\cdot 2^{\binom{k-1}{2}}\,n^{k-1}\,,
\]
 and thus we obtain
 \begin{equation}\label{Eq:A_bounds}
   \sum_{j=1}^k A_j\ \leq\  \tfrac{e^2+1}{2} \, k\,2^{\binom{k-1}{2}}\,n^{k-1}\,.
 \end{equation}

 Setting $l=k{-}2$ and $i=j{-}2$ gives 
\[
  \sum_{j=2}^k B_j\ =\ \sum_{i=0}^{l} 2^{\binom{l-i}{2}}n^{l-i}\tbinom{n+l+1}{i}\ .
\]
 Arguing as in Lemma~\ref{L:technical} bounds this sum by 
\[
   2^{\binom{l}{2}}n^l\   \sum_{i=0}^{\infty} \tfrac{2^{i-1}}{i!}\ =\ 
   \tfrac{e^2}{2} \cdot 2^{\binom{k-2}{2}}\,n^{k-2}\ .
\]
 
 Finally, as $2\leq k\leq n$, we have $n^k-(n-1)^k \leq kn^{n-1}-1$,
which gives 
\[
  \tfrac{1}{2}+2^{\binom{k}{2}} (n^k-(n{-}1)^k)\ \leq\ k 2^{\binom{k}{2}} n^{k-1}\,.
\]
Collecting all these bounds gives the result.
\end{proof}

\subsection{Proof of Lemma~\ref{L:Boundscomp} (1)}
By Lemma~\ref{L:Gen_Gamma}, we have

 \begin{equation}\label{Eq:Gen_Gammabis}
   F_{k-j}\ =\ \biggl(\sum_{\calI\in\binom{[n+k]}{k}^{2^j}}
     \frac{A(\calI) B(\calI)}{p(\calI)}\biggr)
     \cdot \biggl(\,\prod_{i=1}^{n+k}p_i(y)\biggr)^{2^j}\ ,
 \end{equation}
where $A(\calI), B(\calI)$ are real numbers, $\binom{[n+k]}{k}^{2^j}$
consists of lists $\calI$ of $2^j$ elements of $\binom{[n+k]}{k}$, and
$p(\calI)$ is the product over $I$ in $\calI$ of $p_I$, where
if $I=i_1< \cdots <i_k$ then $p_I:=p_{i_1}\cdots p_{i_k}$.

Fix a list $\calI$ and let 
$r$ be the number of times that $p_1$ appears in $p(\calI)$.
Then $0 \leq r \leq 2^j$ and the sparsity of $p_2,\dotsc,p_{n+k}$
implies that the degree of each monomial of 
 \[
    \frac{1}{p(\calI)} \cdot \biggl(\prod_{i=1}^m p_i(y)\biggr)^{2^j}
 \]
lies in the interval
$[(n-1)2^j+r \, , \, n2^j]$.
Lemma~\ref{L:Boundscomp}(1) follows as any value of $r$ in between $0$ and
$2^j$ can occur in a summand of~\eqref{Eq:Gen_Gammabis}.\hfill\qed

\subsection{Proof of Lemma~\ref{L:Boundscomp}(2) and (3)}

Bernstein's Theorem generalizes Kouchnirenko's Theorem to the case when the
polynomials have possibly different supports.
We will only need a special case of Bernstein's Theorem.
The \Blue{{\it Newton polytope}} $\New(f)$ of a polynomial $f$ with integer
exponents is the convex hull of its support. 

\begin{thm}[Bernstein~\cite{Be75}]
 Let $P$ be a fixed polytope in $\R^k$ and 
 suppose that $f_1,\dotsc,f_k$ are polynomials such that for each
 $i=1,\dotsc,k$, the Newton polytope of $f_i$ is contained in 
 $\lambda_i P$, where $\lambda_i\in\N$.
 Then the number of non-degenerate solutions in $(\C^\times)^k$ to the system 
 $f_1(y)=\cdots=f_k(y)=0$ is at most 
 $\lambda_1\dotsb\lambda_k\cdot k! \vol(P)$.
\end{thm}

%
We prove Lemma~\ref{L:Boundscomp}(3).
By Lemma~\ref{L:Boundscomp}(1), if we set
 \[
   P\ :=\ \{(i_1,i_2\dotsc,i_k) \in \R_>^k \mid n-1 \leq i_1+i_2+\cdots+i_k \leq
   n\} \,
 \]
then the Newton polytope of $F_{j}$ is contained in $2^{k-j} P$ for
$j=1,\dotsc,k$.
By Bernstein's Theorem the number of non-degenerate solutions to the system 
$F_1=\cdots=F_k=0$ in $\Delta \subset (\R^\times)^k\subset (\C^\times)^k$ 
is at most
 \[
    \left(\prod_{j=1}^{k} 2^{k-j}\right) k! \vol(P)\ =\ 
     2^{\binom{k}{2}} k! \vol(P)\,.
 \]
To finish the proof of Lemma~\ref{L:Boundscomp}(3), we note that 
$\vol(P)=\frac{1}{k!}(n^k-(n-1)^k)$.
\qed
\medskip

We now prove Lemma~\ref{L:Boundscomp}(2),
adapting the proof of Lemma~\ref{L:Bounds}(2) by
taking into account of the sparsity of the polynomials $F_1,\dotsc,F_k$. 
Define $\mu_j$ as in~\eqref{Eq:mu_def} by the equations
$f_1(y)=\dotsb=f_{j-1}(y)=1$,
where the functions $f_i$ are defined by~\eqref{Eq:Comp_Gale}.
Then the statement of Lemma~\ref{L:mu} holds if we now claim that
$\overline{\mu_j}\cap\partial\overline{\Delta}$ is a union of $(k{-}j)$-faces of
$\overline{\Delta}$ which forms a PL-submanifold of $\overline{\Delta}$, except
possibly at the origin.
As before, $2\flat(C_j)=|\overline{C_j}\cap\partial\overline{\Delta}|$.
Since $\overline{C_j}$ is the subset of $\overline{\mu_j}$ where we also have
 \begin{equation}\label{E:algebraic}
   \overline{F}_{j+1}(y)\ =\ \overline{F}_{j+2}(y)\ =\ \cdots\ =\
   \overline{F}_k(y)\ =\ 0\,,
 \end{equation}
we estimate $|\overline{C_j}\cap\partial\overline{\Delta}|$ by estimating the
number of points $y$ in the union of $(k{-}j)$-faces of $\overline{\Delta}$
where~\eqref{E:algebraic} holds.

First, $\mu_k=C_k$, so $2\flat(C_k)\leq \Phi_0(\overline{\Delta})$.
Otherwise, suppose that $1\leq j\leq k{-}1$ and let $\phi$ be a 
$(k{-}j)$-face of $\overline{\Delta}$.

If $\phi$ is a linear face, then 
we may apply a linear transformation and assume that the facets supporting
$\phi$ are defined by $y_i=0$ for $i>k{-}j$, so that the face $\phi$ affinely
spans the coordinate subspace $\R^{k-j}$ with coordinates $y_1,\dotsc,y_{k{-}j}$.
Thus the points of $\overline{C}_j\cap \phi$ are a subset of
the points of~\eqref{E:algebraic} in $\R^{k-j}$ having non-zero
coordinates. 
This homogeneous linear transformation does not alter the sparsity of the
polynomials $F_i$ or $\overline{F_i}$.

By Lemma~\ref{L:Boundscomp}(1), the substitution of $y_i=0$ for $i>k{-}j$ into
  $\overline{F}_{k-l}$ gives a polynomial with Newton polytope a subset of
  $2^{l}P$, where  
 \[
   P\ :=\ \{(i_1,\dotsc,i_{k-j}) \in \R_>^{k-j} \mid
    n-1 \leq i_1+i_2+\cdots+i_{k-j} \leq n\} \,.
 \]
By Bernstein's Theorem, the number of  points of $\overline{C}_j\cap \phi$ is
at most 
 \[
    2^{\binom{k-j}{2}}  (k{-}1)! \cdot \vol(P)\ =\ 
    2^{\binom{k-j}{2}}\left(n^{k-j}-(n-1)^{k-j} \right)\,.
\]

If $\phi$ is not a linear face, then we simply apply B\'ezout's Theorem to bound
the number of solutions to the system~\eqref{E:algebraic} lying in $\phi$ and
obtain 
 \begin{equation}\label{E:solfornotlinear}
  \biggl(\prod_{i=1}^{k-j} n2^{k-j-i} \biggr)\ =\ 2^{\binom{k-j}{2}}n^{k-j}\ .
 \end{equation}
This completes the proof of Lemma~\ref{L:Boundscomp}(2).\qed

\subsection{Concrete bounds when $k=2, 3$}
When $k=2$, $\overline{\Delta}$ is a polygon in $\R^2$
and the estimate~\eqref{Eq:Kp_form} becomes
\[
   2\kappa(f)\ \leq\ |V(\psi_1,\psi_2)| \ \leq\ \flat(C_2)\ +\ \flat(C_1)\ +\ 
   |V(\Gamma_2,\Gamma_1)|\,.
\]
Lemma~\ref{L:Boundscomp} provides the bounds
 \[
  \begin{array}{lll}
   2\flat(C_2) & \leq & \Phi_0(\overline{\Delta})\,,\\\rule{0pt}{13pt}
   2\flat(C_1) & \leq &
   \Phi_1^\ell(\overline{\Delta})+n\Phi_1^{n\ell}(\overline{\Delta})\,,
     \qquad\mbox{and} \\
   |V(\Gamma_2,\Gamma_1)| & \leq & 4n-2\,.\rule{0pt}{13pt}
  \end{array}
 \]

\begin{thm} 
 Let $f$ be a polynomial with $n{+}3$ monomials in $n$ variables with $2\leq n$.
 Then $\overline{\Delta}$ is a polygon with
 $\Phi_1^\ell(\overline{\Delta}) \leq 2$, $\Phi_1^{n\ell}(\overline{\Delta}) \leq 2$
 and $\Phi_0(\overline{\Delta})\leq 4$.
 Thus $\kappa(f) \leq \lfloor\frac{5n+1}{2}\rfloor$.
\end{thm}

When $k=3$, $\overline{\Delta}$ is a polytope in $\R^3$
and the estimate~\eqref{Eq:Kp_form} becomes
\[
   2\kappa(f)\ \leq\ 
   |V(\psi_1,\psi_2,\psi_3) | \ \leq\ \flat(C_3) + \flat(C_2)\ +\ \flat(C_1)\ +\ 
   |V(\Gamma_3, \Gamma_2,\Gamma_1)|\,.
\]
Lemma~\ref{L:Boundscomp} provides the bounds
 \begin{eqnarray*}
  2\flat(C_3) & \leq & \Phi_0(\overline{\Delta})\,,\\
  2\flat(C_2) & \leq &\Phi_1^\ell(\overline{\Delta}) +
  n\Phi_1^{n\ell}(\overline{\Delta})\,,\\ 
  2\flat(C_1) & \leq & 2 (2n-1)
  \Phi_2^\ell(\overline{\Delta})+2n^2\Phi_2^{n\ell}(\overline{\Delta})\,,
  \qquad\mbox{and}\\  
  |V(\Gamma_3, \Gamma_2,\Gamma_1)| & \leq & 24n^2-24n+8\,.
 \end{eqnarray*}
%

\begin{thm}\label{T:boundscomponentsthree}
 Let $f$ be a polynomial with $n+4$ monomials and $n\geq 2$ variables.
Then $\overline{\Delta}$ is a $3$-polytope with 
\begin{itemize}

\item $\Phi_0(\overline{\Delta})  \leq 4n+4$

\item $\Phi_1^\ell(\overline{\Delta}) \leq n+2$, \ 
    $\Phi_1^{n\ell}(\overline{\Delta}) \leq 2n+5$,

\item $\Phi_2^\ell(\overline{\Delta}) \leq n+2$, \ 
      $\Phi_2^{n\ell}(\overline{\Delta}) \leq 2$.
\end{itemize}
Thus 
 \[
     \kappa(f)\ \leq\ \tfrac{29}{2}n^2 - 8n + \tfrac{9}{2}\,.
 \]
\end{thm}

\begin{proof}
We need only to establish the bounds on the number of faces.
Let $K$ be the cone in $\R^3$ defined by $p_i(y)\geq 0$ for $2\leq i\leq n+3$.
It is combinatorially equivalent to the cone over the face $\phi^\infty$ at
infinity, which is a polygon with at most $n{+}2$ edges.
Thus 
\[
  \Phi_0(\overline{K})\ \leq\ n+3\,,\qquad
  \Phi_1^\ell(\overline{K})\,,\; 
  \Phi_1^{n\ell}(\overline{K})\,,\; 
  \Phi_2^\ell(\overline{K})\ \leq\ n+2\,,\qquad\mbox{and}\qquad
  \Phi_2^{n\ell}(\overline{K})\ \leq\ 1\,.
\]

We cut $K$ with the half space $p_1(y)>0$ to obtain $\Delta$.
This can only decrease the number of linear faces, so we consider the 
other faces.
This adds at most one facet, so $\Phi^{n\ell}_2(\overline{\Delta})\leq 2$.
We can assume that $\Delta\neq K$ and let $\phi^1$ be the new facet.

Since $\phi^1$ cannot meet more than $2n{+}2$ of the edges of $K$, and must cut
off at least one vertex of $K$, we have $\Phi_0(\overline{\Delta})\leq 4n+4$.

Since $K$ has $n+3$ facets, $\Phi^{n\ell}_1(\overline{\Delta})\leq 
\Phi_2^\ell(\overline{K})+n+3$, which completes the proof.
\end{proof}

\providecommand{\bysame}{\leavevmode\hbox to3em{\hrulefill}\thinspace}
\providecommand{\MR}{\relax\ifhmode\unskip\space\fi MR }
\providecommand{\MRhref}[2]{%
  \href{http://www.ams.org/mathscinet-getitem?mr=#1}{#2}
}
\providecommand{\href}[2]{#2}

\end{document}